\documentclass[12pt]{article}
\usepackage{amsmath,amssymb,amsthm,eucal,mathrsfs}
\usepackage[pdftex]{graphicx}
\usepackage{hyperref}
\usepackage{color}
\usepackage[normalem]{ulem}
\usepackage{enumerate}

\font\tenrm=cmr10

\font\bigss=cmssdc10 scaled 2300

\font\cmsslll=cmss10 at 14 pt

\renewcommand{\a}{\alpha}
\renewcommand{\b}{\beta}

\renewcommand{\d}{\delta}

\renewcommand{\k}{\kappa}

\renewcommand{\o}{\omega}

\newcommand{\bC}{\mathbb{C}}
\newcommand{\bR}{\mathbb{R}}

\renewcommand{\square}{\kern1pt\vbox
               {\hrule height 0.6pt\hbox{\vrule width 0.6pt\hskip 3pt
    \vbox{\vskip 6pt}\hskip 3pt\vrule width 0.6pt}\hrule height0.6pt}
    \kern1pt}

\newcommand{\ra}{\rightarrow}

\newcommand{\dif}{\mathrm d}

\newtheorem{Pb}{Problem}

\newtheorem{Th}{Theorem}
\newtheorem{Prop}[Th]{Proposition}

\newtheorem{Cor}[Th]{Corollary}
\newtheorem{Lem}[Th]{Lemma}
\newtheorem{Def}[Th]{Definition}

\theoremstyle{definition}

\newcommand{\bP}{\begin{Pb}\ \ }
\newcommand{\eP}{\end{Pb}}
\newcommand{\bt}{\begin{Th}\ \ }
\newcommand{\et}{\end{Th}}
\newcommand{\bp}{\begin{Prop}\ \ }
\newcommand{\ep}{\end{Prop}}
\newcommand{\bc}{\begin{Cor}\ \ }
\newcommand{\ec}{\end{Cor}}
\newcommand{\bl}{\begin{Lem}\ \ }
\newcommand{\el}{\end{Lem}}
\newcommand{\bd}{\begin{Def}\ \ }
\newcommand{\ed}{\end{Def}}
\newcommand{\pf}{\begin{proof}[{\it Proof:\ \ }]}
\newcommand{\epf}{\end{proof}}

\newcommand{\be}{\begin{equation}}
\newcommand{\ee}{\end{equation}}

\newcommand\re[1]{(\ref{#1})}
\newcommand{\arr}{\begin{array}{rlll}}
\newcommand{\ea}{\end{array}}
\newcommand{\bea}{\begin{eqnarray}}
\newcommand{\eea}{\end{eqnarray}}
\newcommand{\bean}{\begin{eqnarray*}}
\newcommand{\eean}{\end{eqnarray*}}

\catcode`@=11
\@addtoreset{equation}{section}
\catcode`@=12

\theoremstyle{remark}

\newtheorem{Rem}{Remark}  
\newcommand{\br}{\begin{Rem}\ \ }  
\newcommand{\er}{\end{Rem}}   
\begin{document}
\rightline{}
\vskip 1.5 true cm
\begin{center}
{\bigss  Quarter-pinched Einstein metrics 
interpolating\\[.5em] between real and complex hyperbolic metrics}
\vskip 1.0 true cm
{\cmsslll  V.\ Cort\'es and A.\ Saha
} \\[3pt]
{\tenrm   Department of Mathematics\\
and Center for Mathematical Physics\\
University of Hamburg\\
Bundesstra{\ss}e 55,
D-20146 Hamburg, Germany\\
vicente.cortes@uni-hamburg.de, arpan.saha@uni-hamburg.de}\\[1em]
\vspace{2ex}
{November 14, 2017}
\end{center}
\vskip 1.0 true cm
\baselineskip=18pt
\begin{abstract}We show that the one-loop quantum deformation 
of the universal hypermultiplet provides a family of complete
$1/4$-pinched negatively curved quaternionic K\"ahler (i.e.\ half conformally flat Einstein) metrics 
$g^c$, $c\ge 0$, on $\bR^4$.  The metric $g^0$ is the complex hyperbolic metric 
whereas the family $(g^c)_{c>0}$ is equivalent to a family of metrics $(h^b)_{b>0}$ depending
on $b=1/c$  and smoothly extending to $b=0$ for which $h^0$ is the real hyperbolic metric.  
In this sense the one-loop deformation interpolates between the real and the complex hyperbolic metrics. 
We also determine the (singular) conformal structure at infinity for the above families. 
\\[.5em] 
\end{abstract}

\section*{Introduction}\label{sec:Intro}
Einstein deformations of rank one symmetric spaces of non-compact type  have been considered
by various authors, see \cite{P,L,B1,B2} and references therein. In particular, LeBrun has shown that the quaternionic hyperbolic metric 
on the smooth manifold $\bR^{4n}$ admits deformations by complete quaternionic K\"ahler metrics. These metrics
are constructed using deformations of the twistor data and depend on functional parameters. However, the sectional curvature of the deformed metrics does not seem to have been studied. 

In previous work \cite{ACM,ACDM} a geometric construction of a 
class of quaternionic K\"ahler manifolds of negative scalar curvature was described. The manifolds 
in this class are obtained from projective special K\"ahler manifolds and come in 
one-parameter families. In string theory, such families can be interpreted as perturbative quantum corrections 
to the hypermultiplet moduli space metric \cite{RSV}.  
The one-parameter families are known as one-loop deformations of the supergravity c-map 
metrics. The simplest example corresponds to the case when the initial projective special K\"ahler manifold is a point. 
In that case one obtains the family of metrics  
\begin{equation}\label{eq:1ldUHmetric}
    \begin{split}
        g^c = {1\over4\rho^2}&\left[\frac{\rho + 2c}{\rho + c}\,\mathrm d\rho^2 + \frac{\rho + c}{\rho + 2c}(\mathrm d \tilde\phi + \zeta^0\mathrm d\tilde\zeta_0 - \tilde\zeta_0\mathrm d\zeta^0)^2\right.\\ 
&\quad \left. + 2(\rho + 2c)\left((\mathrm d\tilde\zeta_0)^2 + (\mathrm d\zeta^0)^2\right)\right],
    \end{split}
\end{equation}
where  $(\rho, \tilde{\phi}, \zeta^0, \tilde\zeta_0)$ are standard coordinates on the manifold $M:=\bR^{>0}\times \bR^3\cong \bR^4$ and $c\ge 0$. 
This is a deformation of the complex hyperbolic metric $g^0$ (known as the 
universal hypermultiplet metric in the physics literature \cite{RSV}) by complete quaternionic K\"ahler\footnote{Recall that in dimension four  quaternionic K\"ahler manifolds are defined as half conformally flat Einstein manifolds.} metrics, see \cite[Remark 8]{ACDM}. 
Using the c-map and its one-loop deformation it is also possible to deform higher rank quaternionic K\"ahler symmetric spaces 
and, more generally, quaternionic K\"ahler homogeneous spaces by families of complete quaternionic K\"ahler metrics 
depending on one or several parameters \cite{CDS,CDJL}.

In this paper we prove that the metrics \re{eq:1ldUHmetric} are all negatively curved and $\frac14$-pinched, see 
Theorem \ref{mainthm}. 
By similar calculations, we also show 
that Pedersen's deformation of the real hyperbolic $4$-space\footnote{This deformation is induced 
by a deformation of the 
standard conformal structure of $S^3$ at the boundary of the real hyperbolic space by a rescaling along the fibres
of the Hopf fibration \cite{P}.}, which depends on a parameter $m^2\ge 0$, has
negative 
curvature  if $m^2 <1$, see Theorem \ref{PedersenmetricThm}.
These are presumably the first examples of non-locally symmetric 
complete Einstein four-manifolds of negative curvature. For the family \re{eq:1ldUHmetric}, we show  in Section \ref{cinftySec} that the limit 
$c\ra \infty$ is well-defined
after a suitable change of coordinates and parameter, and that it is given by the real hyperbolic metric. Furthermore,
we perform another change of coordinates in order to analyze the conformal structure at infinity. We find  in Section \ref{confSec} that the 
conformal structure induced by $g^c$ (for $0< c<\infty$) on the boundary sphere $S^3$ is \emph{singular}  
precisely at a single point $p_\infty$, which we can consider as the south pole, where it has a double pole. The point 
$p_\infty$ is also a special point with respect to the asymptotic behaviour of the metric. In fact, the metric 
$g^c$ (considered as a metric on the $4$-ball $B^4$ with boundary $S^3$) is asymptotic to the real hyperbolic 
metric on the complement in $B^4$ of any neighborhood of $p_\infty$ but it is not 
near $p_\infty$. These observations show that the family of metrics $g^c$ cannot be obtained as an Einstein deformation
induced by a deformation of the conformal structure at the boundary in the spirit of \cite{B1}. 

\subsubsection*{Acknowledgements} We are very grateful to Alexander Haupt for checking some of our calculations.  Furthermore, we thank 
Olivier Biquard, Gerhard Knieper and Norbert Peyerimhoff for helpful comments. 
This work was supported by the German Science Foundation (DFG) under the Research Training Group 1670 ``Mathematics inspired by String Theory". 
Finally, the authors would like to express a special thanks to the  Mainz Institute for 
Theoretical Physics (MITP) 
for its hospitality and support. 

\section{The limit $c\ra \infty$}\label{cinftySec}
We introduce a second one-parameter family of metrics given by
\begin{equation}
	\begin{split}
		h^b = {1\over4\rho'^2}&\left[\frac{b\rho' + 2}{b\rho' + 1}\,\mathrm d\rho'^2 + \frac{b\rho' + 1}{b\rho' + 2}(\mathrm d \tilde\phi' + b\zeta'^0\mathrm d\tilde\zeta'_0 - b\tilde\zeta'_0\mathrm d\zeta'^0)^2\right.\\ 
&\quad \left. + 2(b\rho' + 2)\left((\mathrm d\tilde\zeta'_0)^2 + (\mathrm d\zeta'^0)^2\right)\right],
	\end{split}
\end{equation} 
where $b>0$. This is in fact equivalent to the one-loop deformation $g^c$ for $c>0$ under the identifications $c=1/b$ and $(\rho,\tilde\phi,\zeta^0,\tilde\zeta_0) = (\rho',\tilde\phi',\sqrt b\,\zeta'^0,\sqrt b\,\tilde\zeta'_0)$. But now the family can be extended to the $b=0$ case. This implies that after the above parameter-dependent coordinate transformation the $c\rightarrow\infty$ limit of the one-loop deformation $g^c$ is indeed well-defined and is given by the metric 
\begin{equation}
	\begin{split}
		h^0 = {1\over4\rho'^2}&\left[2\,\mathrm d\rho'^2 + \frac{1}{2}\,\mathrm d \tilde\phi'^2 + 4(\mathrm d\tilde\zeta'_0)^2 + 4(\mathrm d\zeta'^0)^2\right], 
	\end{split}
\end{equation}
which has constant curvature $-2$.
\section{Asymptotics and conformal structure at infinity}\label{confSec}
We would like to determine the conformal structure of the family of metrics $(g^c)_{c\ge 0}$ on the sphere at the boundary of $M$. In our coordinates, this consists of the hyperplane at $\rho = 0$, along with a point at infinity  $p_\infty$. In order to be able to directly see the singularity at $p_\infty$, we consider the following change of coordinates: 
\begin{equation}\label{eq:changeofcoord1}
\begin{split}
\rho &= \Re\left(\frac{1-z_1}{1+z_1}\right) - \left\lvert\frac{z_2}{z_1 + 1}\right\rvert^2 = \frac{1-|z_1|^2 -|z_2|^2}{|z_1+1|^2}
,\\
\tilde \phi &= -\Im \left(\frac{1-z_1}{1+z_1}\right),\quad \zeta := \zeta^0 + \mathrm i \tilde\zeta_0 = \frac{\sqrt{2}\,z_2}{1+z_1}.
\end{split}
\end{equation}
This is indeed a diffeomorphism from $M=\mathbb R_{>0} \times  \mathbb \bR^3=\mathbb R_{>0} \times \mathbb R\times \mathbb C$ to the unit ball $B_{\mathbb C}^2$ in $\bC^2$, as it admits the following (smooth) inverse:
\begin{equation}
\begin{split}
z_1 &= \frac{1-\left(\rho + |\zeta|^2/2 -\mathrm i\tilde\phi\right)}{1+\left(\rho + |\zeta|^2/2 -\mathrm i\tilde\phi\right)},\quad
z_2 = \frac{\sqrt{2}\,\zeta}{1+\left(\rho + |\zeta|^2/2 -\mathrm i\tilde\phi\right)}.
\end{split}
\end{equation}
As a result of the above change of coordinates, the boundary is mapped to the unit sphere $S^3\subset \mathbb C^2$, and $p_\infty$ is mapped to the south pole $(z_1,z_2) = (-1,0)$. 
\bp
In the coordinates introduced in \eqref{eq:changeofcoord1}, the conformal structure at the boundary $[g^c|_{\partial M}]$, for $c>0$ is singular at $p_\infty$ ($z_1 = -1$) and away from the singularity is given by the nondegenerate conformal structure:
\begin{equation}
\begin{split}
\left[g^{c}|_{\partial M}\right] &=
\left[\left(2\, \Re\left(\mathrm d\left(\frac{1-z_1}{1+z_1}\right) - \left(\frac{2\,\overline z_2}{1+ \overline z_1}\right)\mathrm d\left(\frac{z_2}{1+z_1}\right) \right)^2\right.\right.\\
 &\quad +\frac{1}{2}\, \Im\left(\mathrm d\left(\frac{1-z_1}{1+z_1}\right) - \left(\frac{2\,\overline z_2}{1+ \overline z_1}\right)\mathrm d\left(\frac{z_2}{1+z_1}\right) \right)^2\\
&\quad +\left.\left.\left.  8c\left\lvert\mathrm d \left(\frac{z_2}{1+z_1}\right)\right\rvert^2\right)\right\rvert_{\partial M}\right].
\end{split}
\end{equation}
Meanwhile the conformal structure for $c=0$ is supported only on the CR distribution $\mathscr D$ on $S^3$ and is given by
\begin{equation}
\left[g^{0}|_{{\mathscr D}\times {\mathscr D}}\right] = \left[\left.\left(\left\lvert\mathrm d \left(\frac{z_2}{1+z_1}\right)\right\rvert^2\right)\right\rvert_{{\mathscr D}\times {\mathscr D}}\right].
\end{equation}
\ep
\begin{proof}
For any $c\ge 0$, the metric $g^c$ in the new coordinates is given by
\begin{equation}\label{eq:metricnewcoord}
\begin{split}
g^c &= \frac{1}{4\rho^2}\left[\frac{\rho + 2c}{\rho + c}\, \Re\left(\mathrm d\left(\frac{1-z_1}{1+z_1}\right) - \left(\frac{2\,\overline z_2}{1+ \overline z_1}\right)\mathrm d\left(\frac{z_2}{1+z_1}\right) \right)^2\right.\\
&\quad\left. +\frac{\rho + c}{\rho + 2c}\, \Im\left(\mathrm d\left(\frac{1-z_1}{1+z_1}\right) - \left(\frac{2\,\overline z_2}{1+ \overline z_1}\right)\mathrm d\left(\frac{z_2}{1+z_1}\right) \right)^2+ 4(\rho + 2c)\left\lvert\mathrm d \left(\frac{z_2}{1+z_1}\right)\right\rvert^2\right],
\end{split}
\end{equation}
where now $\rho = \frac{1-|z_1|^2 -|z_2|^2}{|z_1+1|^2}$ is considered as a function of $(z_1,z_2)$. 
The above metric is well-defined and nondegenerate when $|z_1|^2 +|z_2|^2 < 1$. Moreover we see that for $c>0$, the conformal structure at the boundary $[g^c|_{\partial M}]=[(4\rho^2g^c)|_{\partial M}]$ is singular at $z_1 = -1$. Away from the singularity, it may be computed to be the following:
\begin{equation*}
\begin{split}
\left[g^{c}|_{\partial M}\right] &=
\left[\left(2\, \Re\left(\mathrm d\left(\frac{1-z_1}{1+z_1}\right) - \left(\frac{2\,\overline z_2}{1+ \overline z_1}\right)\mathrm d\left(\frac{z_2}{1+z_1}\right) \right)^2\right.\right.\\
 &\quad +\left.\left.\left.\frac{1}{2}\, \Im\left(\mathrm d\left(\frac{1-z_1}{1+z_1}\right) - \left(\frac{2\,\overline z_2}{1+ \overline z_1}\right)\mathrm d\left(\frac{z_2}{1+z_1}\right) \right)^2 +  8c\left\lvert\mathrm d \left(\frac{z_2}{1+z_1}\right)\right\rvert^2\right)\right\rvert_{\partial M}\right].
\end{split}
\end{equation*}
Meanwhile, in the case $c=0$, the (rescaled) metric in \eqref{eq:metricnewcoord} becomes:
\begin{equation}
\begin{split}\label{decompEq}
\rho g^0 &= \frac{1}{4\rho}\left\lvert\mathrm d\left(\frac{1-z_1}{1+z_1}\right) - \left(\frac{2\,\overline z_2}{1+ \overline z_1}\right)\mathrm d\left(\frac{z_2}{1+z_1}\right) \right\rvert^2+ \left\lvert\mathrm d \left(\frac{z_2}{1+z_1}\right)\right\rvert^2.
\end{split}
\end{equation}
The second term stays finite at the boundary but the first term blows up, except on its kernel, which may be verified to be spanned by the following two vector fields:
\begin{equation*}
\overline z_2\, \frac{\partial}{\partial z_1} - \left(\frac{1 - |z_2|^2 +\overline z_1}{1+z_1}\right)\frac{\partial}{\partial z_2},\quad z_2\, \frac{\partial}{\partial \overline z_1} - \left(\frac{1 - |z_2|^2 + z_1}{1+ \overline z_1}\right)\frac{\partial}{\partial \overline z_2}.
\end{equation*}
At the boundary, the above become vector fields spanning the CR distribution $\mathscr D$  on $S^3$:
\begin{equation*}
\overline z_2\, \frac{\partial}{\partial z_1} - \overline  z_1\,\frac{\partial}{\partial z_2},\quad z_2\, \frac{\partial}{\partial \overline z_1} - z_1\,\frac{\partial}{\partial \overline z_2}.
\end{equation*}
The conformal structure at the boundary  $\left[g^{0}|_{{\mathscr D}\times {\mathscr D}}\right]$ 
is defined as the nondegenerate conformal structure on ${\mathscr D}$ obtained by keeping only the finite term
in the above decomposition \re{decompEq}, see \cite{B1}. 
Thus the conformal structure $\left[g^{0}|_{\partial M}\right]$ is supported only on the CR distribution ${\mathscr D}$  and is given by
\begin{equation*}
\left[g^{0}|_{{\mathscr D}\times {\mathscr D}}\right] := \left[(\rho g^{0})|_{{\mathscr D}\times {\mathscr D}}\right] = \left[\left.\left(\left\lvert\mathrm d \left(\frac{z_2}{1+z_1}\right)\right\rvert^2\right)\right\rvert_{{\mathscr D}\times {\mathscr D}}\right].
\end{equation*}
\end{proof}
We would also like to determine the conformal structure of the family of metrics $(h^b)_{b\ge 0}$ on the sphere at the boundary of $M$. As in the case of $g^c$ above, in order to directly see the singularity at the point at infinity $p_\infty$, we again carry out a change of coordinates that maps $M=\mathbb R_{>0} \times \mathbb R^3$ to the unit ball $B_{\mathbb R}^4$ in $\bR^4$:
\begin{equation}\label{eq:changeofcoord2}
    \begin{split}
        \rho' &= {1- w^2 - x^2 - y^2 -z^2\over (1+ w)^2 + x^2 + y^2 + z^2},\quad
        \tilde\phi' = {4x\over (1+ w)^2 + x^2 + y^2 + z^2},\\
        \zeta'^0 &= {\sqrt{2}\,y\over (1+ w)^2 + x^2 + y^2 + z^2},\quad
        \tilde\zeta_0' = {\sqrt{2}\,z\over (1+ w)^2 + x^2 + y^2 + z^2}.
    \end{split}
\end{equation}
This is indeed a diffeomorphism, with (smooth) inverse given by 
\begin{equation}
    \begin{split}
        w &= {1- \rho'^2 - \tilde\phi'^2/4  -2\left(\zeta'^0\right)^2 - 2\,\tilde\zeta_0'^2\over (1+ \rho')^2 + \tilde\phi'^2/4 +2\left(\zeta'^0\right)^2+ 2\,\tilde\zeta_0'^2},\quad
        x = {\tilde\phi'\over (1+ \rho')^2 + \tilde\phi'^2/4 +2\left(\zeta'^0\right)^2+ 2\,\tilde\zeta_0'^2},\\
        y &= {2\sqrt{2}\,\zeta'^0\over  (1+ \rho')^2 + \tilde\phi'^2/4 +2\left(\zeta'^0\right)^2+ 2\,\tilde\zeta_0'^2},\quad
        z = {2\sqrt{2}\,\tilde\zeta_0'\over  (1+ \rho')^2 + \tilde\phi'^2/4 +2\left(\zeta'^0\right)^2+ 2\,\tilde\zeta_0'^2}.\\
    \end{split}
\end{equation}
As a result of the above change of coordinates, the boundary is mapped to the unit sphere $S^3\subset \mathbb R^4$, and the point at infinity $p_\infty$ is mapped to the south pole $(w,x,y,z) = (-1,0,0,0)$. 

\bp
In the coordinates introduced in \eqref{eq:changeofcoord2}, the conformal structure $[h^b|_{\partial M}]$ on the boundary sphere for $b>0$ is singular at $w=-1$ and away from the singularity is given by
\begin{equation}
\begin{split}
\left[h^b|_{\partial M}\right]
&= \left[\left(\frac{1}{2}\left(\mathrm d \left(\frac{2x}{1+w}\right) + \frac{b}{2}\left(\frac{y}{1+w}\right)\mathrm d\left(\frac{z}{1+w}\right) - \frac{b}{2}\left(\frac{z}{1+w}\right)\mathrm d\left(\frac{y}{1+w}\right)\right)^2\right.\right.\\ 
&\quad \left.\left.\left. + 2\left(\left(\mathrm d\left(\frac{y}{1+w}\right)\right)^2 + \left(\mathrm d\left(\frac{z}{1+w}\right)\right)^2\right)\right)\right\rvert_{\partial M}\right].
    \end{split}
\end{equation}
Moreover, for $b=0$, the conformal structure $\left[h^0|_{\partial M}\right]$ is the standard conformal structure on $S^3$.
\ep
\begin{proof}
At the boundary and away from the south pole, we have $w^2 + x^2 + y^2 + z^2 = 1$ and $w\neq -1$. So, the restrictions of the coordinate functions $\rho',\tilde{\phi}',\zeta'^0,\tilde\zeta'_0$ to $\partial M$ are given as functions of $w,x,y,z$ as follows:
\begin{equation}
\rho'|_{\partial M} = 0,\quad \tilde{\phi}'|_{\partial M} = \frac{2x}{1+w}, \quad \zeta'^0|_{\partial M} = \frac{y}{\sqrt 2\,(1+w)}, \quad \tilde\zeta'_0|_{\partial M} = \frac{z}{\sqrt 2\,(1+w)}.
\end{equation}
A straightforward substitution therefore yields
\begin{equation}
\begin{split}
(4\rho'^2h^b)|_{\partial M} 
 &= \left(\frac{1}{2}\left(\mathrm d \left(\frac{2x}{1+w}\right) + \frac{b}{2}\left(\frac{y}{1+w}\right)\mathrm d\left(\frac{z}{1+w}\right) - \frac{b}{2}\left(\frac{z}{1+w}\right)\mathrm d\left(\frac{y}{1+w}\right)\right)^2\right.\\ 
&\quad +\left.\left. 2\left(\left(\mathrm d\left(\frac{y}{1+w}\right)\right)^2 + \left(\mathrm d\left(\frac{z}{1+w}\right)\right)^2\right)\right)\right\rvert_{\partial M}.
\end{split}
\end{equation}
The conformal structure is nondegenerate with a double pole at the south pole (i.e.~$w=-1$) for $b>0$. When $b=0$, the above becomes:
\begin{equation}
    \begin{split}
        \left[h^0|_{\partial M}\right] &=\left[\left.\left(2\left(\mathrm d\left(\frac{x}{1+w}\right)\right)^2 +  2\left(\mathrm d\left(\frac{y}{1+w}\right)\right)^2 + 2\left(\mathrm d\left(\frac{z}{1+w}\right)\right)^2\right)\right\rvert_{\partial M}\right]\\
        &= \left[\left.\left(\frac{2(\dif x^2 + \dif y^2 + \dif z^2)}{(1+w)^2}+\frac{2(x^2 + y^2 + z^2)\,\dif w^2 }{(1+w)^4}-\frac{4(x\,\dif x+y\,\dif y + z\,\dif z)\,\dif w}{(1+w)^3}\right)\right\rvert_{\partial M}\right]\\
        &= \left[\left.\left(\frac{2(\dif x^2 + \dif y^2 + \dif z^2)}{(1+w)^2}+\frac{2(1- w^2)\,\dif w^2 }{(1+w)^4}+\frac{4(w\,\dif w)\,\dif w}{(1+w)^3}\right)\right\rvert_{\partial M}\right]\\
        &= \left[\left.\left(\frac{2(\dif w^2+\dif x^2 + \dif y^2 + \dif z^2)}{(1+w)^2}\right)\right\rvert_{\partial M}\right]= \left[\left.\left(\dif w^2+\dif x^2 + \dif y^2 + \dif z^2\right)\right\rvert_{\partial M}\right].
    \end{split}
\end{equation}
This is the standard conformal structure on $S^3$ i.e.~the conformal class to which the restriction of the Euclidean metric on $\mathbb R^4$ to $S^3$ belongs.
\end{proof}

\section{Computation of the curvature tensor}
Our goal is to prove the following result (restated in Theorem \ref{mainthmrep}):
\bt \label{mainthm}  For the 
one-loop deformation $g^c$, $c>0$, 
the pinching function $p\mapsto \d_p$ defined in \re{pinching} satisfies $\frac{1}{4} < \d < 1$  and attains the 
boundary values asymptotically when  $\tilde{\rho} =\rho/c$ approaches $0$ or $\infty$, respectively, which is to say, $M$ is everywhere 
(at least) ``quarter-pinched".
\et
In order to prove this, we first compute the curvature associated with the metric in \eqref{eq:1ldUHmetric} by making use of the Cartan formalism. 
In this formalism, we choose an orthonormal frame $(e_I)_{I=1,\ldots ,4}$ and denote the dual co-frame by $(\theta^I)$ so that $g^c=\sum_{I}\theta^I\otimes \theta^I$. The way we have presented the metric in \eqref{eq:1ldUHmetric} suggests an obvious choice, namely
\begin{equation}\label{eq:theta}
    \begin{split}
        \theta^1&:=F(\rho)\,\mathrm d\rho,\quad\,\,\, \theta^2:=G(\rho)(\mathrm d\tilde{\phi}+\zeta^0\mathrm d\tilde\zeta_0-\tilde\zeta_0\mathrm d\zeta^0),\\
        \theta^3&:=H(\rho)\,\mathrm d\tilde\zeta_0,\quad
        \theta^4:=H(\rho)\,\mathrm d\zeta^0.
    \end{split}
\end{equation}
where $F(\rho), G(\rho), H(\rho)$ are functions of $\rho$ given by
\begin{equation}\label{eq:F}
\begin{split}
F(\rho)&=\frac{1}{2\rho}\,\sqrt{\frac{\rho + 2c}{\rho + c}},\quad
G(\rho)=\frac{1}{2\rho}\,\sqrt{\frac{\rho + c}{\rho + 2c}},\quad
H(\rho)=\frac{\sqrt{2(\rho + 2c)}}{2\rho}.
\end{split}
\end{equation}
The $\mathfrak{so}(4)$-valued connection $1$-form $\omega = (\omega^I_J)$ and curvature $2$-form $\Omega = (\Omega^I_J)$ corresponding to the Levi-Civita connection $\nabla$ and its curvature tensor $R$ 
are defined by 
\begin{equation}\label{eq:relation}
        \nabla_{v}e_I = \sum_J\omega^J_I(v)e_J, \quad \Omega^J_I(v,w)=g^c(R(v,w)e_I,e_J), 
\end{equation}
for any vector vector fields $v$ and $w$. The forms $\omega^J_I$ and $\Omega^J_I$ can be calculated through the Cartan structural equations: 
\begin{equation}\label{eq:Cartan}
    \begin{split}
        \mathrm d\theta^I &= \sum_J\theta^J\wedge \omega^I_J,\quad
        \mathrm d\omega^I_J = \Omega^I_J + \sum_K \omega^K_J\wedge\omega^I_K.\\
    \end{split}
\end{equation}
In fact, the first equation is equivalent to the vanishing of torsion and determines the forms $\omega^J_I=-\o^I_J$ uniquely. 
We now gather together the results of the calculation in the following two lemmata. We omit the proofs, which consist of just checking
the structure equations.
\bl
The connection $1$-forms $\omega^I_J$ in \eqref{eq:Cartan} are given by
\begin{equation}\label{eq:omega}
    \begin{split}
        \omega^1_2 &= -\omega^2_1 = \frac{1}{F(\rho)}{2\rho^2 + 5c\rho +4c^2\over 2\rho(\rho+c)(\rho+2c)}\,\theta^2,\quad
        \omega^1_3 = -\omega^3_1 = \frac{1}{F(\rho)}{\rho +4c\over 2\rho(\rho+2c)}\,\theta^3,\\
        \omega^1_4 &= -\omega^4_1 = \frac{1}{F(\rho)}{\rho +4c\over 2\rho(\rho+2c)}\,\theta^4,\quad\quad\quad\quad
        \omega^2_3 = -\omega^3_2 = -\frac{1}{F(\rho)}{1\over 2(\rho +2c)}\,\theta^4,\\
        \omega^2_4 &= -\omega^4_2 = \frac{1}{F(\rho)}{1\over 2(\rho +2c)}\,\theta^3,\,\,\,\quad\quad\quad\quad
        \omega^3_4 = -\omega^4_3 = \frac{1}{F(\rho)}{1\over 2(\rho +2c)}\,\theta^2.
    \end{split}
\end{equation}
\el
\bl\label{lem:Omega}
The curvature $2$-forms $\Omega^I_J$ in \eqref{eq:Cartan} are given by 
\begin{equation}\label{eq:Omega}
    \begin{split}
        \Omega^1_2 &= -\Omega^2_1 = -A_{\mathrm{I}}(\rho)\,\theta^1\wedge\theta^2+2A_{\mathrm{III}}(\rho)\,\theta^3\wedge\theta^4,\\
        \Omega^1_3 &= -\Omega^3_1 = -A_{\mathrm{II}}(\rho)\,\theta^1\wedge\theta^3+A_{\mathrm{III}}(\rho)\,\theta^2\wedge\theta^4,\\
        \Omega^1_4 &= -\Omega^4_1 =-A_{\mathrm{II}}(\rho)\theta^1\wedge\theta^4-A_{\mathrm{III}}(\rho)\,\theta^2\wedge\theta^3,\\
        \Omega^2_3 &= -\Omega^3_2 = -A_{\mathrm{III}}(\rho)\,\theta^1\wedge\theta^4-A_{\mathrm{II}}(\rho)\,\theta^2\wedge\theta^3,\\
        \Omega^2_4 &= -\Omega^4_2 =~~\,A_{\mathrm{III}}(\rho)\,\theta^1\wedge\theta^3-A_{\mathrm{II}}(\rho)\,\theta^2\wedge\theta^4,\\
        \Omega^3_4 &= -\Omega^4_3 =\,\,2A_{\mathrm{III}}(\rho)\,\theta^1\wedge\theta^2-A_{\mathrm{I}}(\rho)\,\theta^3\wedge \theta^4,
    \end{split}
\end{equation}
where $A_{\mathrm{I}}$, $A_{\mathrm{II}}$, and $A_{\mathrm{III}}$ are given by
\begin{equation}
    \begin{split}
        A_{\mathrm{I}}(\rho) &:= \frac{4\rho^3 + 12c\rho^2 + 24c^2\rho + 16c^3}{(\rho + 2c)^3},\\
        A_{\mathrm{II}}(\rho)  &:= \frac{\rho^3 + 12c\rho^2 + 24c^2\rho + 16c^3}{(\rho + 2c)^3},\\
        A_{\mathrm{III}}(\rho)  &:= -\frac{\rho^3}{(\rho + 2c)^3}.
    \end{split}
\end{equation}
\el

\section{Eigenspaces of the curvature operator}

In this section we consider the curvature operator $\mathscr{R}: \Lambda^2TM \rightarrow \Lambda^2TM$ which is defined by
\[ g^c(\mathscr{R}X\wedge Y, Z\wedge W) = g^c(R(X,Y)W,Z),\]
where on the left-hand side $g^c$ denotes the scalar product on bi-vectors which is induced by the Riemannian metric $g^c$:
\[ g^c(X\wedge Y, Z\wedge W) = g^c(X,Z)g^c(Y,W)-g^c(X,W)g^c(Y,Z).\]
Identifying vector with co-vectors by means of the metric, we will consider the curvature operator as a map 
\be \label{curvopEq}\mathscr{R}: \Lambda^2T^*M \rightarrow \Lambda^2T^*M.\ee 
As such it maps  $\theta^I\wedge \theta^J$ to 
$\Omega^I_J$. 
The endomorphism $\mathscr R$ is self-adjoint with respect to (the metric on $\Lambda^2T^*M$ induced by) $g^c$. It follows, that all
eigenvalues are real and that there exists an orthonormal eigenbasis.  
\bp \label{specProp}The following (anti-)self-dual $2$-forms 
\begin{equation}
        \alpha^\pm_{JKL} = \theta^1\wedge\theta^J \pm \theta^K\wedge\theta^L, 
\end{equation}
where $(J,K,L)$ is a cyclic permutation of $(2,3,4)$, 
form an eigenbasis of the curvature operator \re{curvopEq}  of the one-loop deformation \eqref{eq:1ldUHmetric}. 
The corresponding eigenvalues $\lambda^\pm_{JKL}$ are 
\begin{equation}
    \begin{split}
        \lambda_{234}^+ &= -2\left[1 + 2\left(\frac{\rho}{\rho +2c}\right)^3\right],\\
            \lambda_{234}^- = \lambda_{342}^- = \lambda_{423}^- &= -2,\\
            \lambda_{342}^+ = \lambda_{423}^+ &= -2\left[1 - \left(\frac{\rho}{\rho +2c}\right)^3\right].
    \end{split}
\end{equation}
In particular, when $c\neq 0$, the above depends only on the ratio $\tilde \rho:= \rho/c$: 
\begin{equation}
    \begin{split}
       \lambda_{234}^+ &= -2\left[1 + 2\left(\frac{\tilde\rho}{\tilde\rho +2}\right)^3\right],\\
           \lambda_{234}^- = \lambda_{342}^- = \lambda_{423}^- &= -2,\\
           \lambda_{342}^+ = \lambda_{423}^+ &= -2\left[1 - \left(\frac{\tilde\rho}{\tilde\rho +2}\right)^3\right].
    \end{split}
\end{equation}
\ep

\pf
From Lemma \ref{lem:Omega} we see that $\mathscr R$ is block diagonal, whereby the bundle $\Lambda^2T^\ast M$ of $2$-forms decomposes into three invariant subbundles $\Lambda^2_{234}T^* M$, $\Lambda^2_{342}T^* M$, and $\Lambda^2_{423}T^* M$, where $\Lambda^2_{JKL}T^* M$ denotes the span of $\theta^1\wedge\theta^J$ and $\theta^K\wedge\theta^L$.  By inspection, we may read off the two eigen-$2$-forms $\alpha^\pm_{JKL}$ in $\Lambda^2_{JKL}T^\star M$.  The corresponding eigenvalues are
\begin{equation}
    \begin{split}
        \lambda_{234}^+ 
        = -A_{\mathrm{I}} + 2A_{\mathrm{III}} = -\frac{6\rho^3 + 12c\rho^2 + 24c^2\rho + 16c^3}{(\rho + 2c)^3}&=-2\left[1 + 2\left(\frac{\rho}{\rho +2c}\right)^3\right],\\
        \lambda_{234}^- 
        =-A_{\mathrm I} - 2A_{\mathrm{III}} = -\frac{2\rho^3 + 12c\rho^2 + 24c^2\rho + 16c^3}{(\rho + 2c)^3}&=-2,\\
        \lambda_{342}^- = \lambda_{423}^- 
        =-A_{\mathrm{II}} + A_{\mathrm{III}} = -\frac{2\rho^3 + 12c\rho^2 + 24c^2\rho + 16c^3}{(\rho + 2c)^3}&=-2,\\
        \lambda_{342}^+ = \lambda_{423}^+ 
        =-A_{\mathrm{II}} - A_{\mathrm{III}} = -\frac{12c\rho^2 + 24c^2\rho + 16c^3}{(\rho + 2c)^3}&=-2\left[1 -\left(\frac{\rho}{\rho +2c}\right)^3\right].
    \end{split}
\end{equation}
\epf
From the above computation, we may read off the Ricci curvature $\mathrm{Rc}:T^*M\rightarrow T^*M$ and Weyl curvature $\mathscr W:\Lambda^2T^*M\rightarrow\Lambda^2T^*M$ as follows:
\begin{equation}
	\begin{split}
		\mathrm{Rc}&:=\sum_{I=1}^4 \iota(e_I)\circ\mathscr R\circ\varepsilon(\theta^I)=-6\,\mathrm{id}_{T^*M},\\
		\mathscr W &:= \mathscr R -\frac{1}{2}\, \mathrm{Rc}\wedge\mathrm{id}_{T^*M}+\frac{1}{3}\,\mathrm{tr}(\mathscr R)\,\mathrm{id}_{\Lambda^2T^*M}\\
		 &=\left(\mathscr R + 2\,\mathrm{id}_{\Lambda^2T^*M}\right)=\frac{1}{2}\,(1+\star)\left(\mathscr R + 2\,\mathrm{id}_{\Lambda^2T^*M}\right),
	\end{split}
\end{equation}
where $\iota(e^I)$ and $\varepsilon(\theta^I)$ are the interior and exterior products respectively, $\star$ is the Hodge star operator, and $\mathrm{Rc}\wedge\mathrm{id}_{T^*M}:\Lambda^2T^*M\rightarrow\Lambda^2T^*M$ is an endomorphism given by 
\begin{equation*}
(\mathrm{Rc}\wedge\mathrm{id}_{T^*M})(\theta^I\wedge \theta^J)=\mathrm{Rc}(\theta^I)\wedge\theta^J - \mathrm{Rc}(\theta^J)\wedge\theta^I.
\end{equation*}
Thus, we see that the metric $g^c$ is Einstein and its Weyl curvature $\mathscr W$ is self-dual, that is, $(M,g^c)$ is indeed quaternionic K\"ahler.
\section{Sectional curvature and pinching  of the one-loop deformation}
Since any element of $\Lambda^2TM$ can be written as a linear combination
of eigenvectors of $\mathscr R$, the sectional curvature 
\begin{equation*}
 K(\Pi) =  g^c(\mathscr{R}u\wedge v, u\wedge v)
\end{equation*}
of a plane $\Pi \subset TM$ with orthonormal basis $(u,v)$ 
can be written as a convex linear combination of the eigenvalues of $\mathscr R$. So the spectrum of $\mathscr R$, determined
in Lemma \ref{specProp}, shall provide bounds on $K$.

In order to obtain the pointwise maximum and minimum of the sectional curvature one has to minimise and maximise 
$g^c(\mathscr{R} \a , \a)$   subject to the conditions $\alpha\wedge\alpha = 0$ (decomposability) and $g^c(\a , \a ) =1$.
This leads us to the following lemma
\bl\label{lem:maxmin}
For any point $p\in M$, we have the following bounds for the sectional curvature of the 
one-loop deformation  \eqref{eq:1ldUHmetric}{\rm :}
\begin{equation}
    \begin{split}
        \max_{\Pi \subset T_pM} K(\Pi) &= {1\over 2}(\max\{\lambda^+_{234}(p),\lambda^+_{342}(p),\lambda^+_{423}(p)\}\\
        &\quad + \max\{\lambda^-_{234}(p),\lambda^-_{342}(p),\lambda^-_{423}(p)\}),\\
        \min_{\Pi \subset T_pM} K(\Pi) &= {1\over 2}(\min\{\lambda^+_{234}(p),\lambda^+_{342}(p),\lambda^+_{423}(p)\}\\
        &\quad + \min\{\lambda^-_{234}(p),\lambda^-_{342}(p),\lambda^-_{423}(p)\}).
    \end{split}
\end{equation}
\el
\begin{proof}
We consider a general $2$-form $\alpha$ written in terms of the eigen-$2$-forms as follows
\begin{equation}\label{eq:alpha}
        \alpha = \sum_{\epsilon, (J,K,L)}a^\epsilon_{JKL}\alpha^\epsilon_{JKL},
\end{equation} 
where $(J, K, L)$ runs over the cyclic permutations of $(2,3,4)$, and $\epsilon$ runs over the values $\pm$. By decomposing $\a$ into its self-dual and anti-self-dual parts, we see that two equations $\alpha\wedge\alpha = 0$ and $g^c(\a , \a ) =1$ are together equivalent to
\begin{equation}\label{constraintEq}
    \begin{split}
        (a^+_{234})^2 + (a^+_{342})^2 + (a^+_{423})^2 = \frac{1}{4},\\
        (a^-_{234})^2 + (a^-_{342})^2 + (a^-_{423})^2 = \frac{1}{4}.
    \end{split}
\end{equation}
On plugging \eqref{eq:alpha} into $g^c(\mathscr{R} \a , \a)$ , we find that
\begin{equation} 
    \begin{split}
        K(\Pi) &= \frac12 \left[4(a^+_{234})^2\lambda^+_{234} + 4(a^+_{342})^2\lambda^+_{324} + 4(a^+_{423)})^2\lambda^+_{423}\right]\\
        &\quad + \frac12\left[4(a^-_{234})^2\lambda^-_{234} + 4(a^-_{342})^2\lambda^-_{324} + 4(a^-_{423})^2 \lambda^-_{423}\right].
    \end{split}
\end{equation}
Under the constraint \re{constraintEq}, the expressions within the square brackets are each convex combinations of three eigenvalues of $\mathscr R$. Therefore in order to maximise or minimise $K(\Pi)$ we need to respectively maximise or minimise these convex combinations separately. 
\end{proof}
In the limit $\tilde \rho\rightarrow 0$, all the eigenvalues become $-2$ as for the real hyperbolic space $\bR\mathbf H^4$ with constant negative sectional curvature $-2$. Meanwhile, in the limit $\tilde \rho\rightarrow \infty$, the pointwise maximum of the sectional curvature is $-1$ and the pointwise minimum is $-4$, giving a pinching of $1/4$ as for the complex hyperbolic plane $\bC\mathbf H^2$.

The interpolation of the pinching between these two limits is described in the following proposition. 
\bp
The pointwise pinching of the metric $g^c$ for $c>0$ at a point $p=(c\tilde\rho,\tilde\phi,\tilde\zeta_0,\zeta^0)\in M$ is given by
\begin{equation} \label{pinching}
\d_p := {\max \{ K(\Pi) \mid  \Pi \subset T_pM\} \over \min \{ K(\Pi) \mid  \Pi \subset T_pM\} }
 =\frac{\tilde\rho^3 +12\tilde\rho^2 + 24 \tilde\rho + 16 }{4\tilde\rho^3 +12\tilde\rho^2 + 24 \tilde\rho + 16 }.
 \end{equation}
\ep
\begin{proof}
We note that we have $\lambda^+_{234} < \lambda^-_{234} = \lambda^-_{342} = \lambda^-_{423} < \lambda^+_{342)} = \lambda^+_{423}$ for all $\tilde{\rho}>0$. So, we have for all $p\in M$
\begin{equation*}
\begin{split}
\max\{\lambda^+_{234}(p),\lambda^+_{342}(p),\lambda^+_{423}(p)\} &= \lambda^+_{342}(p) = \lambda^+_{423}(p),\\
\min\{\lambda^+_{234}(p),\lambda^+_{342}(p),\lambda^+_{423}(p)\} &= \lambda^+_{234},\\
\max\{\lambda^-_{234}(p),\lambda^-_{342}(p),\lambda^-_{423}(p)\} &= \lambda^-_{234}(p) = \lambda^-_{342}(p)= \lambda^-_{423}(p),\\
\min\{\lambda^-_{234}(p),\lambda^-_{342}(p),\lambda^-_{423}(p)\} &= \lambda^-_{234}(p) = \lambda^-_{342}(p)= \lambda^-_{423}(p).
\end{split}
\end{equation*}
It now follows from Lemma \ref{lem:maxmin} that the pointwise pinching at $p=(c\tilde\rho,\tilde\phi,\tilde\zeta_0,\zeta^0)$ is given by
\begin{equation*}
\d_p 
= \frac{\lambda^+_{342}(p) + \lambda^-_{234}(p)}{\lambda^+_{234}(p)+\lambda^-_{234}(p)}
 =\frac{\tilde\rho^3 +12\tilde\rho^2 + 24 \tilde\rho + 16 }{4\tilde\rho^3 +12\tilde\rho^2 + 24 \tilde\rho + 16 },
\end{equation*}
as was to be shown.
\end{proof}
Now that we have a concrete expression for the pointwise pinching, we can derive our main result (stated in Theorem \ref{mainthm} and restated below):
\bt \label{mainthmrep} For the 
one-loop deformation $g^c$, $c>0$, 
the pinching function $p\mapsto \d_p$ defined in \re{pinching} satisfies $\frac{1}{4} < \d < 1$  and attains the 
boundary values asymptotically when  $\tilde{\rho} =\rho/c$ approaches $0$ or $\infty$, respectively, which is to say, $M$ is everywhere 
(at least) ``quarter-pinched".
\et
\begin{proof}
For any $\tilde\rho > 0$, we see that
\begin{equation}
       1> \d_p = \frac{1}{4} + \frac{9\tilde\rho^2 + 18\tilde\rho + 12}{4\tilde\rho^3 +12\tilde\rho^2 + 24 \tilde\rho + 16 } > \frac{1}{4},
\end{equation}
and that both boundary values are attained asymptotically. 
\end{proof}

\section{Pedersen metric}\label{sec:Pedersen}
We now consider the Pedersen metric defined on the unit ball $B^4_\bR$ as discussed in \cite{P}:
\begin{equation} \label{Pedersenmetric}
\begin{split}
\k^m = \frac{1}{(1-\varrho^2)^2}\left(\frac{1+m^2\varrho^2}{1+m^2\varrho^4}\,\dif\varrho^2+ \varrho^2 (1+m^2\varrho^2)\,(\sigma_1^2+\sigma_2^2)+\frac{\varrho^2(1+m^2\varrho^4)}{1+m^2\varrho^2}\,\sigma_3^2\right),
\end{split}
\end{equation}
where the boundary is the sphere at $\varrho=1$ and $\sigma_1,\sigma_2,\sigma_3$ are the three left-invariant 1-forms on $S^3$ satisfying $\dif\sigma_i=\sum_{j,k}\varepsilon_{ijk}\sigma_j\wedge\sigma_k$. As in the case of the 1-loop  
deformed universal hypermultiplet metric, there is an obvious choice of an orthonormal  co-frame $(\theta^I)$, given by
\begin{equation}
\begin{split}
\theta^1 &= \frac{\varrho}{(1-\varrho^2)}\sqrt{1+m^2\varrho^2}\,\sigma_1,\quad \theta^2 = \frac{\varrho}{(1-\varrho^2)}\sqrt{1+m^2\varrho^2}\,\sigma_2,\\
\theta^3 &= \frac{\varrho}{(1-\varrho^2)}\sqrt{\frac{1+m^2\varrho^4}{1+m^2\varrho^2}}\,\sigma_3,\quad \theta^4 = \frac{1}{(1-\varrho^2)}\sqrt{\frac{1+m^2\varrho^2}{1+m^2\varrho^4}}\,\dif\varrho.
\end{split}
\end{equation}
The steps in the previous sections for the calculation of  the eigenvalues and an eigenbasis  of the curvature  operator $\mathscr R:\Lambda^2T^*M\rightarrow \Lambda^2T^*M$ may be repeated for the Pedersen metric. We summarize the  results in the next proposition.
\bp
The following (anti-)self-dual 2-forms 
\be
\beta^\pm_{IJK}:=\theta^I\wedge\theta^J \pm \theta^K\wedge\theta^4
\ee
where $(I,J,K)$ is a cyclic permutation of $(1,2,3)$, form an eigenbasis of the curvature operator $\mathscr R$  of the 
Pedersen metric \eqref{Pedersenmetric}. The corresponding eigenvalues $\nu^\pm_{IJK}$ are
\be
\begin{split}
        \nu^+_{123}=\nu^+_{231}=\nu^+_{312}&=-4,\\
        \nu^-_{123}&= -4\left(1-\frac{2m^2\left(1-\varrho^2\right)^3}{{\left(m^{2} \varrho^{2} +  1\right)^3}}\right),\\
        \nu^-_{231}=\nu^-_{312}&= -4\left(1+\frac{m^2\left(1-\varrho^2\right)^3}{{\left(m^{2} \varrho^{2} +  1\right)^3}}\right).
\end{split}
\ee
\ep
In order to obtain the pointwise maximum and minimum of the sectional curvature one has to minimise and maximise 
$\k^m(\mathscr{R} \b , \b)$   subject to the conditions $\b\wedge\b = 0$ (decomposability) and $\k^m(\b , \b ) =1$. Again, this calculation proceeds exactly as earlier and so we just summarise the result in the following proposition.
\bp
The pointwise maximum and pointwise minimum of the sectional curvature of the Pedersen metric is given by
\begin{equation}
\begin{split}
	\max_{\Pi \subset T_pM} K(\Pi) &= -4\left(1-\frac{m^2\left(1-\varrho^2\right)^3}{{\left(m^{2} \varrho^{2} +  1\right)^3}}\right),\\
    \min_{\Pi \subset T_pM} K(\Pi) &= -4\left(1+\frac{m^2\left(1-\varrho^2\right)^3}{{2\left(m^{2} \varrho^{2} +  1\right)^3}}\right).
\end{split}
\end{equation}
\ep
In particular, a straightforward rearrangement shows that the pointwise maximum $\max_{\Pi \subset T_pM} K(\Pi)$ becomes nonnegative when the following condition holds:
\begin{equation}
    \varrho^2 \le \frac{\sqrt[3]{m^2}-1}{m^2 + \sqrt[3]{m^2}}.
\end{equation}
Note that this condition cannot hold if $m^2<1$. As a consequence we have the 
following result. 
\bt \label{PedersenmetricThm} The Pedersen metric \eqref{Pedersenmetric} has negative sectional curvature 
if and only if $m^2<1$. For $m^2>  1$ (respectively $m^2=1$) there are negative as well as positive (respectively zero) sectional curvatures near (respectively at) the origin $\varrho=0$. 
\et 


\begin{thebibliography}{ABCD}

\bibitem[ACDM]{ACDM}  D.V.\ Alekseevsky, V.\ Cort\'es, M.\ Dyckmanns and T.\ Mohaupt, {\it
Quaternionic K\"ahler metrics associated with special K\"ahler manifolds}, J.\ Geom.\ Phys.\ {\bf 92} (2015), 271--287.

\bibitem[ACM]{ACM} D.V.\ Alekseevsky, V.\ Cort\'es and T. Mohaupt, {\it Conification of K\"ahler and hyper-K\"ahler manifolds}, Comm. Math. Phys. {\bf 324} 
(2013), no. 2, 637--655.

\bibitem[B1]{B1} O.\ Biquard, {\it 
Einstein deformations of hyperbolic metrics. Surveys in differential geometry: essays on Einstein manifolds}, 235--246, 
Surv.\ Differ.\ Geom., 6, Int.\ Press, Boston, MA, 1999. 

\bibitem[B2]{B2} O.\ Biquard, {\it M\'etriques d'Einstein asymptotiquement sym\'etriques}, 
Ast\'erisque {\bf 265} (2000), vi+109 pp. 

\bibitem[CDJL]{CDJL} V.\ Cort\'es, M.\ Dyckmanns, M.\ J\"ungling and D.\ Lindemann, {\it A class of cubic hypersurfaces and quaternionic K\"ahler manifolds of co-homogeneity one},  arXiv:1701.7882 [math.DG]. 

\bibitem[CDS]{CDS} V.\ Cort\'es, M.\ Dyckmanns and S.\ Suhr, {\it Completeness of projective special K\"ahler and quaternionic K\"ahler manifolds}, 
in ``Special metrics and group actions in geometry,'' Springer INdAM Series vol.\ 23 (2017), arXiv:1607.07232 [math.DG].

\bibitem[L]{L} C.\ LeBrun, {\it On complete quaternionic-K\"ahler manifolds}, Duke Math.\ J.\ {\bf 63} (1991), no.\ 3, 723--743.

\bibitem[P]{P} H.\ Pedersen, {\it Einstein metrics, spinning top motions and monopoles}, Math.\ Ann.\ {\bf 274} (1986), 35--59. 

\bibitem[RSV]{RSV}
  D.~Robles-Llana, F.~Saueressig and S.~Vandoren,
  {\it String loop corrected hypermultiplet moduli spaces},
  JHEP {\bf 0603} 081 (2006).
  \end{thebibliography}
  \end{document}